 \documentclass[11pt]{amsart}

 \usepackage{amssymb}

 \input xy

 \xyoption{all}

 \theoremstyle{plain}

 \newtheorem{theorem}{Theorem}[section]

 \newtheorem{corollary}[theorem]{Corollary}

 \newtheorem{conjecture}[theorem]{Conjecture}

 \newtheorem{proposition}[theorem]{Proposition}

 \newtheorem{fact}[theorem]{Fact}

 \newtheorem{claim}[theorem]{Claim}

 \newtheorem{question}[theorem]{Question}

 \theoremstyle{definition}

 \newtheorem{definition}[theorem]{Definition}

 \newtheorem{remark}[theorem]{Remark}

 \newtheorem{notation}[theorem]{Notation}

 \theoremstyle{remark}

 \newcommand{\Pm}{\mathcal{P}^{-}}

 \newcommand{\cf}{\operatorname{cf}}

 \newcommand{\Aut}{\operatorname{Aut}}

 \newcommand{\tp}{\operatorname{ga-tp}}

 \newcommand{\gaS}{\operatorname{ga-S}}

 \newcommand{\Hanf}{\operatorname{Hanf}}

 \newcommand{\Min}{\operatorname{Min}}

 \newcommand{\Max}{\operatorname{Max}}

 \newcommand{\lan}{\operatorname{L}}

 \newcommand{\Mod}{\operatorname{Mod}}

 \renewcommand{\phi}{\varphi}

 \newcommand{\Union}{\bigcup}

 \newcommand{\initial}\lessdot

 \newcommand{\K}{\operatorname{\mathcal{K}}}

 \newcommand{\id}{\operatorname{id}}

 \newcommand{\LS}{\operatorname{LS}}

 \def\?{?\vadjust

 {\vbox to 0pt{\vskip-7pt\hbox to 1.1\hsize{\hfill\huge ?!}}}}

\begin{document}

 \title[Categoricity in Tame AECS]{Shelah's Categoricity Conjecture from a
 successor for Tame Abstract Elementary Classes}

 \author{Rami Grossberg}
 \email[Rami Grossberg]{rami@andrew.cmu.edu}
 \address{Department of Mathematics\\
 Carnegie Mellon University\\
 Pittsburgh PA 15213}

 \author{Monica VanDieren}
 \email[Monica VanDieren]{mvd@umich.edu}
 \address{Department of Mathematics\\
 University of Michigan\\
 Ann Arbor MI 48109-1109}


 \date{September 16, 2005}

 \begin{abstract}

 We proveÊ a categoricity transfer theorem for tame abstract elementary
 classes.

 \begin{theorem}\label{step up}
 Suppose that $\K$ is a $\chi$-tame abstract
 elementary class and satisfies the amalgamation and joint embedding
 properties and has arbitrarily large models.Ê Let
 $\lambda\geq\Max\{\chi,\LS(\K)^+\}$.Ê If
 $\K$ is categorical in
 $\lambda$ and $\lambda^+$, then $\K$ is categorical in
 $\lambda^{++}$.
 \end{theorem}

 Combining this theorem with some results from \cite{Sh 394}, we derive a
 form of Shelah's Categoricity Conjecture for tame abstract elementary
 classes:

 \begin{corollary}\label{main cat cor}
 Suppose $\K$ is a $\chi$-tame abstract elementary class satisfying the
 amalgamation and joint embeddingÊ properties.
 ÊLet
 $\mu_0:=\Hanf(\K)$.Ê If $\chi\leq\beth_{(2^{\mu_0})^+}$ and
 $\K$ is categorical in some
 $\lambda^+>\beth_{(2^{\mu_0})^+}$, then $\K$ is categorical in $\mu$ for
 all
 $\mu>\beth_{(2^{\mu_0})^+}$.

 \end{corollary}

 \end{abstract}

 \maketitle


 \bigskip

 \section*{Introduction} \label{s:introductionm}

 Let $\K$ be a class of models all of the same language.Ê $\K$ is said to be
 \emph{categorical in a cardinal
 $\mu$} if and only if all the models from $\K$ of cardinality $\mu$ are
 isomorphic.Ê In 1954 Jerzy \L o\' s \cite{Lo}
 conjectured that for a first-order theory $T$ in a countable language its class
 of models $\Mod(T)$ behaves
 like the class of algebraically closed fields of fixed characteristic.
 Namely, if there exists an uncountable cardinal $\lambda$ such that
 $\Mod(T)$ is categorical in $\lambda$ then $\Mod(T)$ is categorical in
 every uncountable cardinal. Michael Morley \cite{Mo} in 1965 published a
 proof of this conjecture andÊ asked about its generalization for
 first-order theories in uncountable languages.Ê SaharonÊ Shelah \cite{Sh
 31} around 1970 managed to prove the relativized \L o\' s conjecture for
 theories in uncountable languages. The work of Morley and Shelah
 introduced a large number of new concepts andÊ devices into model theory
 and established the field known as \emph{stability theory} or by the name
 of
 \emph{classification theory}.


 Already inÊ the forties and fifties Alfred Tarski and Andrzej Mostowski
 realized that first-order logic is too weak to deal with some of the most
 basic objects of mathematics and proposed logics with greater expressive
 power like
 $L_{\omega_1,\omega}$ and $L_{\omega_1,\omega}(\mathbf Q)$.Ê Examples due to
 Morley and
 Jack Silver in the sixties
 gave the impression that while these logics haveÊ good expressive power, there
 is very little structure there
 and almost any conjecture has a counterexample.

 Some ground breaking work of Shelah from the mid seventies (\cite{Sh 48},
 \cite{Sh 87a} and \cite{Sh 87b})
 together with earlier important work by H. Jerome Keisler
 \cite{Ke2} led ShelahÊ Êto realize thatÊ machinery of stability theory
 could be developed for the context of non-first order theories. Around
 1977 Shelah
 proposed a far reaching conjecture to serve as test problem to progress
 in the field.

 \begin{conjecture}[Shelah's Categoricity Conjecture]\label{cat conj}
 Given a countable language $L$ and $T$Ê a
 theory in
 $L_{\omega_1,\omega}$, if $\Mod(T)$ is categorical in $\lambda$ for some
 $\lambda> \beth_{\omega_1}$,\footnote{$\beth_{\omega_1}$ is the Hanf
 number of this logic.}Ê then
 $\Mod(T)$ is categorical in $\chi$ for every
 $\chi\geq \beth_{\omega_1}$.

 \end{conjecture}

 While $L_{\omega_1,\omega}$ surfaced as a general context for developing
 a non-elementary model theory,
 it is not broad enough to capture some mathematically
 significant examples such as Boris Zilber's class of algebraically closed
 fields with pseudo-exponentiation which is connected with Schanuel's
 conjecture of transcendentalÊ number theory
 \cite{Zi}.Ê Furthermore, by focusing on a specific non-first-order logic
 such as $L_{\omega_1,\omega}$ one may be near-sighted, not envisioning
 the underlying stability theory.
 %
 %
 %
 %
 In \cite{Sh 88} the notion of \emph{Abstract Elementary Class} (see
 Definition \ref{def AEC})Ê was introduced by Saharon Shelah. This is a semantic
 generalization of
 $L_{\omega_1,\omega}({\mathbf Q})$-theories, according to \cite{Sh 702}
 AECs are the most general context to have a reasonable model theory.Ê ÊIn
 Shelah own words: ``I have preferred this context, certainly the widest I
 think has any chance at all.''
 The guiding conjecture in the development of a
 classification theory of non-elementary classes is now a strengthening
 ofÊ Conjecture \ref{cat conj} and appears in the list of open problems
 in \cite{Shc}:

 \begin{conjecture}[Shelah's Categoricity Conjecture]\label{cat conj aec}
 Let
 $\K$ be an abstract elementary class.Ê If
 $\K$ is categorical in some $\lambda>\Hanf(\K)$,
 \footnote{For an explanation of $\Hanf(\K)$ see Remark \ref{hanf
 notation}} then for every
 $\mu\geq \Hanf(\K)$, $\K$ is categorical in $\mu$.
 \end{conjecture}

 Similar to the solutions of \L o\' s conjecture, a categoricity
 transfer theorem for non-first-order logic is expected to provide the
 basic conceptual tools necessary for a stability theory for
 non-first-order logic which then may be applied to answer questions in
 other branches of mathematics.Ê Already Shelah's proof of Conjecture
 \ref{cat conj aec} for excellent classes forshadowed tools developed by
 Zilber to study algebraically closed fields with psuedo-exponentiation.

 In \cite{GrVa}, we wanted to develop basic stability theory for AECs with
amalgamation,
 in our attempt to prove certain technical statements necessary for establishing
stability
 spectrum theorem and existence of Morley sequences we introduce a property we
called
 tameness.Ê Later we realized that a relativized version of tameness
for saturated models
appeared
 implicitly in Shelah's proof of his main theorem from \cite{Sh 394}.

 In \cite{GrVa} we introduce the notion of tameness to abstract elementary
classes as a
 context with enough generality to capture many mathematical examples, but
surprisingly poignant enough
 to be accessible from a model theoretic point of view. Zilber's work on
 algebraically closed fields equipped with elliptic curves is tame.Ê An
interesting non-tame example is Shelah and Bradd Hart's example of an
 abstract elementary class which is categorical in
 $\aleph_k$ for every $k<n$, but fails to be categorical in $2^{\aleph_n}$.

categoricity


 While there are over a thousand published pages devoted to a partial
 solution of Conjecture
 \ref{cat conj}, it remains wide open.
 Here we prove an approximation to
 Conjecture
 \ref{cat conj aec} for tame abstract elementary classes.Ê Below we will
 review the history of work towards Conjecture \ref{cat conj aec} noting
 that our result is the most general approximation of Conjecture \ref{cat
 conj aec}.Ê Our proof is unprecedented in providing an upward categoricity
 transfer theorem without employing compactness machinery via manipulations
 of first-order or infinitary syntax.

 Now let us make explicit some of the concepts we have mentioned:

 \begin{definition}\label{def AEC}
 Let $\K$ be a class of structures
 all in the same similarity type $\lan(\K)$, and
 let $\prec_{\K}$ be a partial order on $\K$.Ê The ordered pair
 $\langle \mathcal{K},\prec_{\mathcal{K}}\rangle$ is an
 \emph{abstract elementary class,Ê AEC} for short
 if and only if
 \begin{enumerate}

 \item [A0] (Closure under isomorphism)
 \begin{enumerate}

 \item
 For every $M\in \mathcal{K}$
 and every $\lan(\K)$-structure $N$ if $M\cong N$ then $N\in \mathcal{K}$.
 \item Let $N_1,N_2\in\K$ andÊ $M_1,M_2\in \K$ such that
 there exist $f_l:N_l\cong M_l$ (for $l=1,2$) satisfying
 $f_1\subseteq f_2$ thenÊ $N_1\prec_{\K}N_2$
 implies thatÊ $M_1\prec_{\K}M_2$.

 \end{enumerate}

 \item [A1]
 For all $M,N\in \mathcal{K}$ if $M\prec_{\mathcal{K}} N$ then $M\subseteq N$.
 \end{enumerate}

 \begin{enumerate}
 \item [A2]
 Let $M,N,M^*$ be $\lan(\mathcal{K})$-structures.
 If $M\subseteq N$, $M\prec_{\mathcal{K}} M^*$ and $N\prec_{\mathcal{K}}
 M^*$ then
 $M\prec_{\mathcal{K}} N$.

 \item [A3]
 (Downward L\"owenheim-Skolem)
 There exists a cardinal\\
 $\LS(\mathcal{K})\geq \aleph_0+|\lan(\mathcal{K})|$
 such that
 for every \\
 $M\in \mathcal{K}$
 and for every $A\subseteq |M|$ there exists $N\in \mathcal{K} $ such that
 $N\prec_{\mathcal{K}} M, \; |N|\supseteq A$ and
 $\|N\|\leq |A|+\LS(\mathcal{K})$.

 \item [A4]
 (Tarski-Vaught Chain)
 \begin{enumerate}

 \item
 For every regular cardinalÊ $\mu$ and every \\
 $N\in \mathcal{K}$ if
 $\{M_i\prec_{\mathcal{K}} N\;:\;i<\mu\}\subseteq \mathcal{K}$
 is $\prec_{\K}$-increasing (i.e. $i<j\Longrightarrow M_i\prec_{\mathcal{K}}
 M_j$) then
 $\Union_{i<\mu}M_i\in \mathcal{K}$ and $\Union_{i<\mu}M_i\prec_{\mathcal{K}} N
 $.

 \item
 For every regular $\mu$,
 ifÊ $\{M_i\;:\;i<\mu\}\subseteq \mathcal{K}$ is
 $\prec_{\K}$-increasing then $\Union_{i<\mu}M_i\in \mathcal{K}$ and
 $M_0\prec_{\mathcal{K}} \Union_{i<\mu}M_i$.
 \end{enumerate}
 \end{enumerate}

 For $M$ and $N\in \K$ a monomorphism $f:M\rightarrow N$ is called an
 \emph{$\K$-embedding} if and only if $f[M]\prec_{\K}N$.Ê Thus, $M\prec_{\K}N$
is
 equivalent to ``$\id_M$ is a $\K$-embedding from $M$ into $N$''.
 \end{definition}

 \begin{remark}\label{hanf notation}
 Ê$Hanf(\K)$ is widely accepted abuse of notation.Ê In the case
 $\K=\Mod(\psi)$ for $\psi\in L_{\omega_1,\omega}$ in a countable language
 $\Hanf(\K)=\beth_{\omega_1}$, in the more general case
 $\Hanf(\K)=\beth_{(2^{2^{\LS(\K)}})^+}$ where $\LS(\K)$ is the
 L\"{o}wenheim-Skolem number of the class (the cardinality that appears in
 a Downward-L\"{o}wenheim-Skolem theorem for the class).Ê See \cite{Gr1}
 for a formal definition.

 \end{remark}

 In recent years there has been much activity in several
 concrete
 generalizations of first order model theory.
 They are
 \begin{itemize}
 \item model theory of Banach spaces (see
 \cite{HeIo},\cite{Io1},\cite{Io2})
 \item homogeneous model theory, formerly
 known as finite diagrams stable in power (see
 \cite{Be}, \cite{BeBu}, \cite{BuLe},Ê \cite{GrLe}, \cite{Hy},
 \cite{HySh},
 \cite{Le}, \cite{Sh3},
 \cite{Sh54})Ê and
 \item compact abstract theories (CATs)
 \cite{BY}.
 \end{itemize}
 In all of these contexts, a categoricity transfer theorem has been
 proved.Ê However, these classes are very specialized.Ê AECs and even
 $L_{\omega_1,\omega}$ have much more expressive power than these
 specializations.Ê For instance, the classes of solvable groups and
 universal locally finite groups are AECs but cannot be captured by any of
 these specializations.
 Furthermore, the contexts itemized above are far too limited to
 handle Zilber's class of algebraically closed fields with
 pseudo-exponentiation and its variations.Ê Each of these contexts turns
to be tame.

 We summarize the non-elementary categoricity transfer results known to
 date.
 ÊKeisler in 1971 solved Conjecture \ref{cat conj} under the
 additional assumption of existence of a sequentially homogeneous model
 \cite{Ke1}.Ê Under these same assumptions Olivier Lessmann provides a
 Baldwin-Lachlan style proof of the categoricity transfer result in
 \cite{Le}\footnote{
 The referee pointed out to us that also Hyttinen in \cite{Hy1} obtained a
similar result.
 ÊWe feel that the reader will be interested to know that Lesmann's result is
part of his
 1998 PhD thesis and the relevant paper was submitted for publication and was
widely
 circulated already in 1997. }.
 Unfortunately this is insufficient for Shelah's conjecture, since
 Marcus and Shelah found an example for an
 $L_{\omega_1,\omega}$-sentence which is categorical in all infinite
 cardinals but does not have a sequentially homogeneous model (see
 \cite{Ma}).

 Extending the work of Keisler, Shelah in 1984
 proved Conjecture \ref{cat conj}
 Êfor
 excellent classesÊ [87a] andÊ [87b].Ê The hardest part is to show
 under the assumption of $2^{\aleph_n}<2^{\aleph_{n+1}}$ for all
 $n<\omega$ that $\mbox{I}(\aleph_{n+1},\psi)<\mu(n)$ (for all
 $n$) implies that a certain class of atomic models of a first order
 theory derived from
 $\psi$ is
 excellent.

 Other attempts to prove Shelah's Categoricity Conjecture involved making
 extra set theoretic assumptions. Michael Makkai and Shelah
 \cite{MaSh} proved Conjecture
 \ref{cat conj aec} under the additional assumption that the categoricity
 cardinal
 $\lambda$ is a successor and the class $\K$ is axiomatizable by a
 $L_{\kappa,\omega}$-theory where $\kappa$ is above a strongly compact
 cardinal andÊ $\lambda>\kappa$.
 Oren Kolman and Shelah \cite{KoSh} began the generalization of
 \cite{MaSh} replacing the hypothesis that $\kappa$ is above a strongly
 compact cardinal with the assumption that $\kappa$ is above a measurable
 cardinal.
 Shelah completed this work in \cite{Sh 472}, but only managed to prove a
 partial downward categoricity transfer theorem.

 \begin{fact}
 Let $\K$ be an AEC axiomatized by a $L_{\kappa,\omega}$-theory with
 $\kappa$ measurable.Ê If $\K$ is categorical in some
 $\lambda^+>\Hanf(\K)$, then $\K$ is categorical in every $\mu$ with
 $\Hanf(\K)<\mu\leq\lambda^+$.
 \end{fact}

 The most general context of AECs considered so far are AECs which satisfy
 the amalgamation property.Ê In \cite{Sh
 394}, Shelah proves aÊ partial going down result for these classes:
 %

 \begin{fact}\cite{Sh 394}\label{statement of Sh394}
 Suppose that $\K$ satisfies the amalgamation and joint embedding
 properties.Ê Let $\mu_0:=\Hanf(\K)$.Ê If $\K$ is categorical in some
 $\lambda^+>\beth_{(2^{\mu_0})^+}$, then
 $\K$ is categorical in every $\mu$ such that
 $\beth_{(2^{\mu_0})^+}<\mu\leq\lambda^+$.

 \end{fact}

 One of the better approximations to Shelah's categoricity conjecture for
 AECs can be derived from a theorem due to Makkai and Shelah
 (\cite{MaSh}):

 \begin{fact}[Makkai and Shelah 1990]
 Ê ÊLet $\K$ be an AEC,Ê {Ê $\kappa$ a strongly compact cardinal}
 such that $\LS(\mathcal{K})<\kappa$.Ê Let $\mu_0:=\beth_{(2^\kappa)^+}$.
 If
 $\K$ is categorical in some
 Ê Ê$\lambda^{+}>\mu_0$ then $\K$ is categorical in every $\mu\geq\mu_0$.
 \end{fact}

 It is easy to see (using the assumption that $\kappa$ is strongly compact)
 that any AEC $\K$ as above has the AP (for models of cardinality
 $\ge\kappa$) and is also tame.

 Our main result can be viewed as replacing the assumption of existence of a
strongly
 compact cardinal in the Makkai and Shelah theorem by tameness and the
amalgamation
 property.Ê As a consequence of Corollary \ref{final corollary}, we get:

 \begin{theorem}
 Ê Let $\K$ be an AEC, $\kappa:=\beth_{(2^{\LS(\K)})^+}$.Ê Denote by
 $\mu_0:=\beth_{(2^\kappa)^+}$.Ê Suppose that $\K_{>\kappa}$ has the
 amalgamation property and is tame.
 If
 $\K$ is categorical in some
 Ê Ê$\lambda^{+}>\mu_0$ then $\K$ is categorical in every $\mu\geq\mu_0$.
 \end{theorem}

 This is the first upward categoricity theorem we know in ZFC for AECs.

 Shelah and Villaveces \cite{ShVi} and \cite{Va} begin to study AECs
 with no maximal models under GCH.Ê ÊThe focus of this work is to
 initially prove that the amalgamation property follows from categoricity.
 After informing Shelah about the results presented in our paper, he sent
 an email indicating thatÊ using methods of \cite{Sh 705} (good frames and
 $\Pm(n)$-diagrams) he has made progress towards a categoricity transfer
 theorem for AECs
 under some extra set theoretic
 assumptions.

 %

 All of the upward categoricity transfer results relied heavily on
 syntax, strong compactness or set theoretic assumptions.Ê Until now, no
 upward categoricity result was known (or even suspected in light of the
 Hart and Shelah example) to hold.

 %
 %
 %
 %
 %
 %

 %
 %
 %
 %

 In this paper we extend \cite{Sh 87a}, \cite{Sh 87b} and \cite{Sh 394} by
 presenting an upward categoricity transfer theorem for AECs that
 satisfy the amalgamation property with some level of tameness (see
 Definition
 \ref{tame defn}).

 We thank Andr\'{e}s Villaveces for organizing the Bogot\'{a} Meeting in
 Model Theory 2003 where we initially presented this work.Ê We are indebted
 to the participants of this meeting especially John Baldwin, Olivier
 Lessmann and Andr\'{e}s Villaveces for removing several errors from a
 preliminary draft of this paper, for making
 suggestions and raising questions involving this paper.



 \bigskip

 \section{Background} \label{s:background}

 We assume that $\K$ is an abstract
 elementary class (AEC) and satisfies the amalgamation and joint embedding
 property.Ê Ultimately, we will use a little less, by only assuming
 amalgamation and joint embedding for models of cardinality $\leq\lambda^+$
 where
 $\lambda^+$ is a categoricity cardinal.Ê But for readability we
 make the more global assumptions from the start.

 We will be using the
 basic machinery of abstract elementary classes including Galois-types
 introduced by Shelah in
 \cite{Sh 88} and
 \cite{Sh 394}.Ê For convenience we refer the reader to
 \cite{Gr1} for the definitions and essential results. This work both extends
and
 generalizes some results from
 \cite{Sh 394}.Ê For the reader unfamiliar with \cite{Sh 394} we have
 included statements and definitions of the material that we will use
 explicitly here.

 Following the notation and terminology from \cite{Gr1}, for a class $\K$
 and a cardinal $\mu$ we let
 $
 \K_\mu:=\{M\in
 \K \;:\; \|M\|=\mu\}.
 $

 Let $\mu\geq \LS(\K)$.Ê We say that $\K$ has the \emph{$\mu$-amalgamation
 property} if and only if for any $M_\ell\in \K_\mu$ (for $\ell\in \{0,1,2\}$)
 such
 that $M_0\prec_{\K} M_1$ and $M_0\prec_{\K}M_2$ there are
 $N\in\K_\mu$ and $\K$-embeddings $f_\ell:M_\ell
 \rightarrow N$ such that $f_\ell\restriction M_0=\id_{M_0}$
 for $\ell =1,2$.
 $\K$ has the \emph{amalgamation property} if and only if $\K$ has the
 $\mu$-amalgamation property for all $\mu\geq \LS(\K)$.

 %

 \begin{definition}

 An AEC $\K$ is \emph{Galois-stable in $\mu$} if and only if for every
 $M\in\K_\mu$, the number of Galois-types over $M$ is $\leq\mu$.
 \end{definition}

 \begin{fact}[Claim 1.7(a) of \cite{Sh 394}]\label{cat implies stable}

 If $\K$ is categorical in $\lambda\geq \LS(\K)$, then $\K$ is
 Galois-stable in all $\mu$ with $\LS(\K)\leq\mu<\lambda$.
 \end{fact}

 A slight, but useful, improvement of Fact \ref{cat implies stable} can be
 derived from an upward stability transfer theorem for tame classes which
 appears in
 \cite{BaKuVa}.

 \begin{corollary}\label{stability up}

 If $\K$ is categorical in $\lambda^+$, then $\K$ is Galois-stable in
 $\mu$ for all $\LS(\K)\leq \mu\leq\lambda^+$.

 \end{corollary}

 \begin{remark}

 We can only guarantee stability as high as $\lambda^+$ since we don't
 know whether or not there is tameness for types over larger models (see
 Fact \ref{cat implies tame}).

 \end{remark}

 Working under the amalgamation property, Galois-stability implies the
 existence of Galois-saturated models in much the same way as stability
 implies the existence of saturated models in first order model theory.
 Here we review some facts about Galois-saturated models.

 \begin{definition}

 Let $\mu>\LS(\K)$.

 \begin{enumerate}

 \item
 $M\in\K$ is said to be \emph{$\mu$-Galois-saturated}\footnote{Ê We must
 make the distinction between saturated in first order logic (which is a
 property of $M$ alone) and Galois-saturated models (which is depends on
 $M$ and
 $\K$).} if and only if for every
 $N\prec_{\K}M$ with
 $N\in\K_{<\mu}$ and every Galois-type $p$ over $N$, we have that $p$ is
 realized in
 $M$.Ê A model $M$ is \emph{Galois-saturated} if and only if
 it is $\|M\|$-Galois-saturated.

 \item $M\in\K$ is said to be \emph{$\mu$-model homogeneous} if and only if for
 every
 $N\prec_{\K}M$ with
 $N\in\K_{<\mu}$ and every $N'\in\K_{\|N\|}$ with $N\prec_{\K}N'$ there
 exists a $\K$-mapping $f:N'\rightarrow M$ with $f\restriction N=\id_N$.
 We write $M$ is \emph{model homogeneous} to mean that $M$ is
 $\|M\|$-model homogeneous.
 \end{enumerate}

 \end{definition}

 The following is a central property of classes with the amalgamation
 property:

 \begin{fact}[From \cite{Sh 576}, see also \cite{Gr1}]\label{sat is mh}
 Suppose that $\K$ satisfies the amalgamation property.Ê Let
 $M\in\K_{>\LS(\K)}$.Ê The following are equivalent

 \begin{enumerate}

 \item $M$ is Galois-saturated.

 \item $M$ is model homogeneous.

 \end{enumerate}

 \end{fact}

 The same proof gives aÊ relativized version of Fact \ref{sat is mh} which
 we will use later in this paper:

 \begin{fact}

 Suppose that $\K$ satisfies the amalgamation property.Ê Let
 $M\in\K_{>\LS(\K)}$.Ê The following are equivalent
 \begin{enumerate}

 \item $M$ is $\mu$-Galois-saturated.

 \item $M$ is $\mu$-model homogeneous.

 \end{enumerate}

 \end{fact}

 The following is a well known result that has its origins in
 Bjarni J\'{o}nsson's
 \cite{Jo} precursor of uniqueness of saturated models:

 \begin{fact}\label{uniq mh}

 Let $\K$ be an AEC.Ê Suppose $\mu>\lambda\geq \LS(\K)$.
 If $\K$ is categorical in $\lambda$ then all model-homogeneous models of
 cardinality $\mu$ are isomorphic.

 \end{fact}

 A local relative to model homogeneity is that of being universal over.

 \begin{definition}
 \begin{enumerate}
 \item \index{universal over!$\kappa$-universal over}
 Let $\kappa$ be a cardinal $\geq \LS(\K)$.
 We say $N$ is \emph{$\kappa$-universal over $M$} if and only if
 for every $M'\in\K_{\kappa}$ with $M\prec_{\K}M'$ there exists
 a $\K$-embedding $g:M'\rightarrow N$ such that
 $g\restriction M=\id_{M}$:

 \[
 \xymatrix{\ar @{} [dr] M'
 \ar[drr]^{g}Ê & \\
 M \ar[u]^{\id} \ar[rr]_{\id}
 && N
 }
 \]

 \item \index{universal over}

 We say $N$ is \emph{universal over $M$} or $N$ is \emph{a universal
 extension of $M$} if and only if
 $N$ is $\|M\|$-universal over $M$.

 \item For $M\in\K_\mu$,Ê $\sigma$ a limit ordinal with
 $\sigma\leq \mu^+$ and $M'\in \K_{\mu|\sigma|}$
 we say that $M'$ is a \emph{$(\mu,\sigma)$-limit over $M$}\index{limit
 model!$(\mu,\sigma)$-limit model over $M$} if and only if there exists a
 $\prec_{\K}$-increasing and continuous sequence of models $\langle M_i\in
 \K_{\mu}\mid i<\sigma
 \rangle$
 such that
 \begin{enumerate}

 \item $M= M_0$,

 \item $M'=\Union_{i<\sigma}M_i$
 and

 \item\label{univ cond in defn} $M_{i+1}$ is universal over $M_i$.

 \end{enumerate}

 \end{enumerate}

 \end{definition}

 Notice that in the case $\sigma<\mu^+$ then $M'$ has cardinality $\mu$
 andÊ in the case $\sigma=\mu^+$ we have that $\|M'\|=\mu^+$.

 Fact \ref{exist universal} guarantees the existence of limit models from
 stability and amalgamation assumptions:

 \begin{fact}[\cite{Sh 600}, see \cite{GrVa} for a
 complete proof]\label{exist universal} If $\K$ satisfies the amalgamation
 property and is Galois-stable in $\mu$, then for every $M\in\K_\mu$,
 there exists $N\in\K_\mu$ such that
 $M\prec_{\K}N$ and $N$ is universal over $M$.Ê Thus for any $M\in \K_\mu$
 and $\alpha\leq\mu^+$ there exists $N\succneqq_{\K}M$ which is
 $(\mu,\alpha)$-limit.

 \end{fact}

 Now we switch gears and recall the concept of tameness from
 \cite{GrVa}.

 \begin{definition}\label{tame defn}
 Let $\chi$ be a cardinal number.
 We say the abstract elementary class $\K$ with the
 amalgamation property is
 \emph{$\chi$-tame}\index{$\chi$-tame} provided that for
 $M\in\K_{> \chi}$,
 $p\neqÊ q\in \gaS(M)$Ê implies the existence of $N\prec_{\K}M$
 of cardinality $\chi$ such thatÊ $ p\restriction N\neq q\restriction N$.

 \end{definition}

 A variant of $\chi$-tameness involves limiting the scope of the models.
 \begin{definition}
 Assume $\chi<\mu$.
 We say that $\K$ is $(\mu,\chi)$-tame if and only if for all $M\in\K_{\mu}$ and
 all
 $p,q\in \gaS(M)$ whenever $p\neq q$, then there exists $N\prec_{\K}M$ of
 cardinality $\chi$ such that $p\restriction N\neq q\restriction N$.

 \end{definition}

 \begin{notation}

 When $\K$ is $(\mu,\chi)$-tame for all $\mu\leq\mu'$ we write $\K$ is
 $(\leq\mu',\chi)$-tame.

 \end{notation}

 \begin{fact}[From Main Claim 2.3 of part II on page 288 \cite{Sh 394}]\label{cat implies tame}
 If $\K$ is categorical in some $\lambda^+>\beth_{(2^{\Hanf(\K)})^+}$, then
 $\K$ is
 $(<\lambda^+,\chi)$-tame for all
 $\chi(\Phi^*)\leq\chi<\lambda^+$.

 \end{fact}

 \begin{remark}

 \begin{enumerate}

 \itemÊ We will be using tameness to prove the existence of rooted minimal types
 or monotonicity of minimal
 types (Proposition \ref{monotonicity of
 minimal types}).Ê ÊUnfortunately, \cite{Sh 394} only gives us
 tameness up to and not including the categoricity cardinal $\lambda^+$.

 \item

 Formally, Shelah proves that when two types over a saturated model of
 cardinality
 $\kappa<\lambda^+$ differ, then there is a submodel of cardinality
 $\chi(\Phi^*)$ over which they differ.Ê However, if we assume
 categoricity in
 $\lambda^+>\beth_{(2^{\Hanf(\K)})^+}$, byÊ Fact \ref{statement of Sh394}
 all models of cardinality $\kappa$ are saturated.

 \item

 $\chi(\Phi^*)$ has a formal definitionÊ in \cite{Sh 394} related to
 Ehrenfeucht-Mostowski constructions.Ê Since we will only use the fact that
 $\chi(\Phi^*)$ lies below
 $\Hanf(\K)$, we will not give its
 formal definition.

 \end{enumerate}

 \end{remark}

 \begin{question}
 Does categoricity in $\lambda>\Hanf(\K)$ imply
 $(\lambda,\chi)$-tameness for some $\chi<\lambda$?
 \end{question}


 \bigskip

 \section{Minimal Types} \label{s:min types}

 The main tool in our constructions will be a minimal type which is a
 variation ofÊ Definition $(*)_4$ of Theorem 9.7 of \cite{Sh 394}.

 \begin{definition}

 Let $M\in\K_\mu$ and $p\in \gaS(M)$ be given.Ê We say
 $p$ is
 \emph{minimal} if and only if $p$ is non-algebraic (no $c\in M$ realizes $p$)
 and for
 every $M'\in\K_\mu$ with $M\prec_{\K}M'$, there is exactly one
 non-algebraic extension of $p$ to $M'$.

 \end{definition}

 The proof of the following proposition draws on our tameness
 assumption.Ê We will be interested in applying the
 monotonicity proposition to $\lambda$ where $\lambda$ is the categoricity
 cardinal.Ê Recall that Shelah's work does not guarantee any level of
 tameness in the categoricity cardinal.

 \begin{proposition}[Monotonicity of Minimal Types]\label{monotonicity of
 minimal types} Suppose $\K$ is $(\lambda,\chi)$-tame for some
 $\lambda\geq\mu\geq\chi$.
 If $p\in\gaS(M)$ is minimal with $M\in\K_\mu$, then for all
 $N\in\K_{\lambda}$ extending $M$ and every $q\in \gaS(N)$ extending $p$,
 if
 $q$ is non-algebraic then $q$ is minimal.
 \end{proposition}

 \begin{proof}

 Suppose for the sake of contradiction that $p$ and $q$ are as in the
 statement of the proposition, with $q$Ê non-algebraic but not minimal.
 Since $q$ is not minimal, there exist distinct non-algebraic extensions
 of $q$, say
 $q',q''\in \gaS(N')$ for some
 $N'\in\K_\lambda$ with $N\prec_{\K}N'$.Ê By tameness, we can find
 $M'\in\K_\mu$
 of cardinality $\mu$ such that $M\prec_{\K}M'\prec_{\K}N'$ and
 $q'\restriction M'\neq q''\restriction M'$.Ê Notice that $q'\restriction
 M'$ and
 $q''\restriction M'$ are both non-algebraic extensions of $p$.Ê This
 contradicts the minimality of $p$.
 \end{proof}

 \begin{fact}[Density of Minimal Types \cite{Sh 394}]\label{exist
 minimal} If $\K$ is Galois-stable in $\mu$, then for every $N\in\K_\mu$
 and every
 $q\in\gaS(N)$, there are $M\in\K_\mu$ and $p\in\gaS(M)$ such that
 $N\preceq_{\K}M$, $q\leq p$ and $p$ is minimal.
 \end{fact}

 To prove the extension property for minimal types, we need a few facts
 about non-splitting in AECs.

 \begin{definition}

 Let $\mu>\LS(\K)$ be a cardinal.Ê For $M\in\K$ and
 $p\in \gaS(M)$, we say that \emph{$p$ $\mu$-splits over
 $N$}\index{$\mu$-splits}\index{Galois-type!$\mu$-splits}
 if and only if
 $N\prec_{\K}M$ and there exist $N_1,N_2\in\K_\mu$ and a
 $\prec_{\K}$-mapping $h:N_1\cong N_2$ such that
 \begin{enumerate}

 \item $N\prec_{\K}N_1,N_2\prec_{\K}M$,
 \item $h(p\restriction N_1)\neq p\restriction N_2$
 and

 \item $h\restriction N= \id_N$.

 \end{enumerate}

 \end{definition}

 \begin{remark}
 Consider the $\K$-mapping $h$ in the definition of $\mu$-splitting.
 Notice that we do not require that there is an extension
 $h'\in\Aut(M)$ of
 $h$.
 \end{remark}

 The existence, uniqueness and extension properties for non-splitting
 types have been studied in \cite{Sh 394}, \cite{ShVi} and \cite{Va}.
 Here we state the formulations of these results which we will be using.

 Existence of non-splitting types:

 \begin{fact}[Claim 3.3 of \cite{Sh 394}]\label{nonsplit thm}
 Assume $\K$ is an abstract elementary class and is Galois-stable in $\mu$.
 For every $M\in\K_{\geq\mu}$ and $p\in \gaS(M)$,
 there exists $N\in\K_\mu$ such that $p$ does not $\mu$-split over $N$.
 \end{fact}

 %
 %
 %
 A consequence of the proof of the uniqueness result, Theorem I.4.15 of
 \cite{Va}, is the following:
 \begin{corollary}\label{non split ext of non alg}
 Let $N,M,M'\in\K_\mu$ be
 such that
 $M'$ is universal over
 $M$ and
 $M$ is a limit model over $N$.
 Suppose that $p\in\gaS(M)$ does not $\mu$-split over $N$ and $p$ is
 non-algebraic.Ê For every $M'\in\K$ extending $M$ of cardinality
 $\mu$, if $q\in \gaS(M')$ is an extension of $p$ and does not $\mu$-split
 over $N$,Ê then $q$ is non-algebraic.
 \end{corollary}
 %
 %
 %
 %
 %
 %

 The version of this
 extension and existence result that we will use is the following which
 came about when John Baldwin removed the cofinality requirement in Lemma
 6.3 of
 \cite{Sh 394} withÊ an argument using Ehrenfeucht-Mostowski models. This
 has allowed us to remove the assumption of
 $LS(\K)=\aleph_0$ from previous drafts of this paper.

 \begin{fact}[Corollary 2 of \cite{Ba2}]\label{Baldwin non-split ext}
 Suppose that $\K$ is categorical
 in some
 $\lambda>\LS(\K)$ and $\K$ has arbitrarily large models.Ê Let $\mu$ be a
 cardinal such that
 $\LS(\K)<\mu$ and let $\sigma$ be a limit ordinal with
 $\LS(\K)<\sigma<\mu^+$.
 Then, for every
 $(\mu,\sigma)$-limit model $M$ and every type $p\in \gaS(M)$, there exists
 $N\precneqq_{\K}M$ of cardinality $\mu$ such that for
 every
 $M'\in\K_{\leq\lambda}$ extending $M$, there exists $q\in \gaS(M')$ an
 extension of $p$ such
 that $q$ does not $\mu$-split over $N$.Ê In particular $p$ does not
 $\mu$-split over $N$.
 \end{fact}

 The last property of non-splitting that we will need is monotonicity:

 \begin{proposition}
 If $M_0\prec_{\K}N\prec_{\K}M$ and $p\in \gaS(M)$ does not $\mu$-split
 over $M_0$, then $p\restriction N$ does not $\mu$-split over $M_0$.
 \end{proposition}

 \begin{proof}
 Immediate from the definitions.

 \end{proof}

 Combining the machinery of non-splitting, we identify the following
 relative to Claim 4.3Ê of
 \cite{Sh 576}:

 \begin{proposition}[Extension Property for Minimal Types]\label{exist min
 ext} Suppose that $\K$ has arbitrarily large models.Ê Let $\K$ be
 categorical in some
 $\lambda>\LS(\K)$ and
 $(\lambda,\chi)$-tame for some $\chi<\lambda$.Ê Let $\mu$ be such that
 $\LS(\K)<\mu$. If $p\in \gaS(M)$ is minimal and $M$ is a
 $(\mu,\sigma)$-limit model for some limit ordinal $\LS(\K)<\sigma<\mu^+$,
 then for every $M'\in\K_{\leq\lambda}$ extending $M$, there is a
 minimal $q\in \gaS(M')$ such that $q$ extends $p$.
 \end{proposition}

 \begin{proof}
 Without loss of generality $M'$ is universal over $M$.
 Let $p\in\gaS(M)$ be minimal.
 Since $M$ isÊ $(\mu,\sigma)$-limit model,Ê using Fact
 \ref{Baldwin non-split ext}, we can find a
 proper submodel
 $N\prec_{\K}M$ of cardinality $\mu$ such that for every
 $M'\in\K_{\leq\lambda}$ there exists $q\in \gaS(M')$ extending $p$ such
 that $q$ does not $\mu$-split over $N$.
 Suppose for the sake of contradiction that $q$ is not minimal.Ê Then
 tamenessÊ and
 Proposition \ref{monotonicity of
 minimal types} tells us that $q$ must be algebraic.Ê Let $a\in M'$
 realize $q$ and $M^a\in\K_\mu$ contain $a$ with
 $M\prec_{\K}M^a\prec_{\K}M'$.
 Then $q\restriction M^a$ is also algebraic.
 However, since $q\restriction M^a$ does not
 $\mu$-split over $N$ and extends $p$, by CorollaryÊ \ref{non
 split ext of non alg} we see that
 $q\restriction M^a$ is not-algebraic.Ê This gives us a contradiction.
 \end{proof}

 We now introduce a strengthening of minimal types which allow us to
 transfer Vaughtian pairs from one cardinality to another in the subsequent
 section.Ê The intuition is that there is a small part of a minimal
 types that controls its minimality.

 \begin{definition}
 Let $M\in\K_\mu$ be given.
 A type $p\in \gaS(M)$ is \emph{rooted minimal} if and only if $p$ is minimal
and
 there is
 $N\prec_{\K}M$ of cardinality $<\mu$ such that $p\restriction N$ is
 minimal.
 We say that $N$ is a \emph{root} of $p$.
 \end{definition}

 \begin{proposition}[Existence of rooted minimal types]
 \label{locally min exist prop}
 Let $\K$ be categorical in some $\lambda>\chi^+$ and
 $(\lambda,\chi)$-tame with $\chi\geq\LS(\K)$.
 ThenÊ for every
 $M'\in\K_{\lambda}$, there exists a rooted minimal
 $q\in\gaS(M')$.

 \end{proposition}
 \begin{proof}

 Notice that categoricity in $\lambda$ implies stability in
 $\mu$ with $\LS(K)<\mu<\lambda$ by Fact \ref{cat implies stable}.
 Choose $M\in\K_\mu$ be some $\K$-substructure of $M'$ with
 $\mu\geq\chi$.Ê Since $\K$ is stable in $\mu$ and categorical in
 $\lambda$, we may take $M$ to be a $(\mu,\sigma)$-limit model for some
 limit ordinal $\sigma$ with $\LS(\K)<\sigma<\mu$.
 Furthermore, by Fact \ref{exist
 minimal} and monotonicity of minimal types, we can choose $M$ suchÊ that
 there is
 a minimal type $p\in \gaS(M)$.

 Then by
 Proposition \ref{exist min ext}, there exists a minimal
 $q\in \gaS(M')$ extending $p$.Ê $q$ is rooted.
 \end{proof}

 \begin{remark}
 In Section \ref{s:cat transfer} we will prove the existence of
 rooted minimal types over models of cardinality $\lambda$ when
 $\cf(\lambda)=\omega$ under the assumption that
 $\LS(\K)=\aleph_0$.
 \end{remark}

 \begin{proposition}\label{monotonicity of rooted}
 Suppose $\K$ is $\chi$-tame.Ê Let $N\in\K_{\geq\chi}$.
 If $p\in\gaS(M)$ is rooted minimal with $N\prec_{\K}M$ a submodel such
 that $p\restriction N$ is minimal, then for every $N'$ with
 $N\prec_{\K}N'\prec_{\K}M$ we have that $p\restriction N'$ is minimal.
 \end{proposition}

 \begin{proof}
 Follows by tameness and monotonicity of minimal types.
 \end{proof}

 \bigskip

 \section{Vaughtian Pairs} \label{s:vp}

 Next we prove a Vaughtian pair transfer theorem for
 rooted minimal types.





 \begin{definition}
 Let $\mu\leq\lambda$.Ê Fix $M\in\K_\mu$Ê and
 $p\in\gaS(M)$ a minimal type. A \emph{$(p,\lambda)$-Vaughtian pair} is a
 pair of models
 $N_0,N_1\in \K_\lambda$ such that
 \begin{enumerate}

 \item $M\prec_{\K}N_0\precneqq_{\K}N_1$
 and

 \item no $c\in N_1\backslash N_0$ realizes $p$.

 \end{enumerate}

 \end{definition}

 \begin{fact}[Claim $(*)_8$ of Theorem 9.7 of \cite{Sh 394}]\label{no
 lambda,lambda vps}
 Assume that $\K$ is categorical in some successor cardinal $\lambda^+$.
 If $\lambda>\LS(\K)$, then for every model $M\in\K_{\leq\lambda}$ and
 every minimal type
 $p\in
 \gaS(M)$, there are no
 $(p,\lambda)$-Vaughtian pairs.

 \end{fact}

 %
 %
 %
 %








 \begin{theorem}\label{vp transfer}
 Fix $\mu>\LS(\K)$.
 Let $p$ be a rooted minimal type over a model $M$ of cardinality $\mu$.
 Fix a root $N\prec_{\K}M$ of cardinality $\kappa$, with
 $p\restriction N$ minimal. If
 $\K$ has a
 $(p,\mu)$-Vaughtian pair, then
 there is a
 $(p\restriction N,\kappa)$-Vaughtian pair.

 \end{theorem}

 \begin{proof}

 Suppose that $(N^0,N^1)$ form a $(p,\mu)$-Vaughtian pair.
 Ê Let $C$ denote the set of all
 realizations of
 $p\restriction N$ inside $N^1$.Ê Fix
 $a\in N^1\backslash N^0$.

 We now construct $\langle
 N^0_i,N^1_i\in\K_\kappa\mid i<\kappa^+\rangle$ satisfying the following:

 \begin{enumerate}

 \item $N^0_0=N$

 \item $N^\ell_i\precneqq_{\K}N^\ell$ for $\ell=0,1$

 \itemÊ the sequences $\langle N^0_i\mid i<\kappa^+\rangle$ and
 $\langle N^1_i\mid i<\kappa^+\rangle$ are both
 $\prec_{\K}$-increasing and continuous

 \item $a\in N^1_i\backslash N^0_i$ and

 \item \label{put c in} $C_i:=C\bigcap N^1_i\subseteq N^0_{i+1}$.


 \end{enumerate}

 The construction follows from
 the following:

 \begin{claim}\label{no new realizations in N1}
 If $d\in N^1$ realizes $p\restriction N^0_0$, then $d\in N^0$.Ê Thus
 $C\subseteq N^0$.

 \end{claim}

 \begin{proof}[Proof of Claim \ref{no new realizations in N1}]
 Suppose that $d\in N^1\backslash N^0$Ê realizes $p\restriction N^0_0$.
 Then $\tp(d/N^0)$ is a non-algebraic extension of $p\restriction N^0_0$.
 Since $p\restriction N^0_0$ is minimal, we have that $\tp(d/M)=p$.Ê Since
 $(N^0,N^1)$ form a $(p,\mu)$-Vaughtian pair, it must be the case that
 $d\in N^0$, contradicting our choice of $d$.

 \end{proof}

 The construction is enough:
 Define
 $$E:=\left\{\begin{array}{c|l}{\delta<\kappa^+}& \delta\text{ is a limit
 ordinal},\\
 &\text{for all }i<\delta\text{ and }x\in N^1_i,\\
 &\text{if there exists }
 j<\kappa^+\text{ such that }x\in C_j, \\ &\text{then there exists }
 j<\delta, \text{ such that }x\in C_j\end{array}\right\}.$$Ê Notice that
 $E$ is a club. (We only use the fact that $E$ is non-empty.)Ê Fix
 $\delta\in E$.

 \begin{claim}\label{force vp claim}
 For every $c\in N^1_\delta\cap C$, we have $c\in N^0_\delta$.

 \end{claim}

 \begin{proof}[Proof of Claim \ref{force vp claim}]
 Since $\langle
 N^1_i \mid i<\kappa^+\rangle$ is continuous, there is $i<\delta$ such
 that $c\in N^1_i$.Ê Thus by the definition of $E$, there is a
 $j<\delta$ with $c\in C_j$.Ê By condition $(\ref{put c in})$ of the
 construction, we would have put $c\in N^0_{j+1}\prec_{\K}N^0_\delta$.

 \end{proof}

 %




 Notice that $N^1_\delta\neq N^0_\delta$ since $a\in N^1_\delta\backslash
 N^0_\delta$.
 Thus
 Claim \ref{force vp claim} allows us to conclude that we have
 constructed a
 $(p\restriction N^0_0,\kappa)$-Vaughtian pair
 $(N^0_\delta,N^1_\delta)$.

 \end{proof}

 \begin{corollary}\label{no vp cor}
 Let $\lambda>\LS(K)$.
 If
 $\K$ is categorical inÊ $\lambda$ and $\lambda^+$ and $p$ is
 rooted minimal type over a model of cardinality $\lambda^+$, then there
 are no
 $(p,\lambda^{+})$-Vaughtian pairs.

 \end{corollary}

 \begin{proof}
 Suppose $(N_0,N_1)$ is a $(p,\lambda^+)$-Vaughtian pair.Ê Then by Theorem
 \ref{vp transfer} and Proposition \ref{monotonicity of rooted}, there is
 a $(p\restriction N,\lambda)$-Vaughtian pair where
 $p\restriction N$ is minimal.Ê Since $\K$ is categorical in $\lambda$,
 $N$ is saturated.Ê This contradicts Fact
 \ref{no lambda,lambda
 vps}.

 \end{proof}

 \begin{corollary}\label{many realizations}
 Let $\lambda>\LS(\K)$.
 If
 $\K$ is categorical in $\lambda$ and $\lambda^+$, then every
 rooted minimal type over a model $N$ of cardinality $\lambda^+$ is
 realized
 $\lambda^{++}$ times in every model of cardinality $\lambda^{++}$
 extending $N$.

 \end{corollary}

 \begin{proof}
 Suppose $M\in\K_{\lambda^{++}}$ realizes $p$ only $\alpha<\lambda^+$
 times.

 Let $A:=\{a_i\mid i<\alpha\}$ be an enumeration of the realizations of $p$
 in $M$.Ê We can find $N_0\in\K_{\lambda^+}$ such that
 $N\Union A\subseteq N_0\prec_{\K}M$.Ê Since $M$ has cardinality
 $\lambda^{++}$, we can
 find $N_1\in\K_{\lambda^+}$Ê such that
 $N_0\precneqq_{\K}N_1\prec_{\K}M$.Ê Then $(N_0,N_1)$ form a $(p,
 \lambda^+)$-Vaughtian pair contradicting
 Corollary \ref{no vp cor}.

 \end{proof}


 \bigskip

 \section{Upward Categoricity Transfer Theorems } \label{s:cat transfer}

 The followingÊ Êshows the strength of the assumption of no Vaughtian
 pairs.Ê In order to prove that a model is saturated, it suffices to check
 that the model realizes one rooted minimal type many times.
 Furthermore, Theorem \ref{many realizations of minimal imply sat} provides
 a new sufficient condition for a model to be universal over a submodel.

 \begin{theorem}\label{many realizations of minimal imply sat}
 Suppose
 $M_0\in\K_\lambda$ and
 $r\in\gaS(M_0)$ is a minimal typeÊ such that $\K$ has no
 $(r,\lambda)$-Vaughtian pairs.

 Let $\alpha$ be an ordinal $<\lambda^+$ such that
 $\alpha=\lambda\cdot\alpha$.Ê ÊSuppose $M\in\K_{\lambda}$ has a
 resolution
 $\langle M_i\in\K_{\lambda}\mid i<\alpha\rangle$Ê such that for
 every
 $i<\alpha$, there is
 $c_i\in M_{i+1}\backslash M_i$ realizing $r$.Ê Then $M$ is saturated over
 $M_0$.Ê Moreover if $\K$ is stable in $\lambda$, then $M$ is a
 $(\lambda,\alpha)$-limit model over $M_0$.
 \end{theorem}

 Notice that Proposition \ref{locally min exist prop} and Corollary \ref{no
 vp cor} guarantee that such $r$ and $M_0$ exist when we assume that $\K$
 is tame and categorical in
 $\lambda$ with $\lambda$ a successor cardinal.Ê We use the letter $r$ to
 represent this type since in our applications of this theorem, $r$ will
 be rooted and we will want to distinguish this type from others.

 Theorem \ref{many realizations of minimal imply
 sat} is similar to Claim 5.6 of \cite{Sh 576}.Ê It is also related to
 a result of
 \cite{ShVi} and
 \cite{Va} that the top of a relatively full tower of length
 $\alpha=\lambda\cdot\alpha$ is a $(\mu,\alpha)$-limit model.

 \begin{proof}[Proof of Theorem \ref{many realizations of minimal imply
 sat}] Let $\langle M_i\mid i<\alpha\rangle$ and $r$ be given as in the
 statement of the theorem. Without loss of generality, we can assume
 that the resolution $\langle M_i\mid i<\alpha\rangle$ is continuous.
 Fix $q\in\gaS(M_0)$.
 We will prove that $M$ realizes $q$
 by
 constructing a
 model $M'$ which realizes $q$ and a $\K$-mapping from $M$ into this model
 $M'$ simultaneously.
 Then we will show that this is in fact an
 isomorphism.

 Since $\alpha=\lambda\cdot\alpha$, we can fix a
 collection of disjoint sets $\{S_i\mid i<\alpha\}$ such that
 $\alpha=\Union_{i<\alpha}S_i$Ê and each $S_i$ is unbounded in $\alpha$
 of cardinality
 $\lambda$ and
 $\Min(S_i)\geq i$.

 We define by induction on $i<\alpha$ sequences of models
 $\langle N'_i\mid i<\alpha\rangle$ and $\langle M'_i\mid i<\alpha\rangle$
 and a sequence of $\K$-mappings $\langle h_i\mid i<\alpha\rangle $.
 Additionally, for $i<\alpha$ we fix a sequence $\langle a_\zeta\mid
 \zeta\in S_i\rangle$.Ê We require:
 \begin{enumerate}
 \item $M'_0$ realizes both $r$ and $q$,

 \item $N'_i\prec_{\K}M'_i$,

 \item $\langle N'_i\mid i<\alpha\rangle$ and $\langle
 M'_i\mid i<\alpha\rangle$ are $\prec_{\K}$-increasing and
 continuous sequences of models in $\K_\lambda$,

 \item $N'_0=M_0$,

 \item $\langle a_\zeta\mid\zeta\in S_i\rangle$ is an enumeration of
 $\{a\in M'_i\mid a\models r\}$,

 \item $a_i\in N'_{i+1}$ ,

 \item $h_i:M_i\cong N'_i$,

 \item $\langle h_i\mid i<\alpha\rangle$ is increasing and continuous with
 $h_0=\id_{M_0}$
 and

 \item when $\K$ is stable in $\lambda$, we additionally require
 $M'_{i+1}$ is universal over
 $M'_i$.

 \end{enumerate}

 The construction is possible:Ê For $i=0$, we take $N'_0=M_0$ and let
 $M'_0$ be a
 extension of $M_0$ of cardinality $\lambda$ realizing $r$ and $q$.Ê If
 possible we choose $M'_0$ to be a
 universal extension of $M_0$ of cardinality $\lambda$.Ê Set
 $h_0=\id_{M_0}$.Ê Let $\langle a_\zeta\mid\zeta\in S_0\rangle$ be some
 (possibly repeating) enumeration of $\{a\in M'_0\mid a\models r\}$.

 Suppose that we have defined for all $k\leq j$, $N'_k,M'_k, h_k$ and
 $\langle a_\zeta\mid\zeta\in S_k\rangle$.Ê Let $a_j$ be given.Ê Notice
 that
 $a_j$ has been definedÊ since $\min S_l> j$ for
 $l\geq i$.Ê If $a_j$ is already in $N'_j$, then we simply amalgamate the
 following diagram
 \[
 \xymatrix{\ar @{} [dr] M_{j+1}
 \ar[r]^{f}Ê &M^{*} \\
 M_j \ar[u]^{\id} \ar[r]_{h_j}
 & M'_{j} \ar[u]_{\id}
 }
 \]
 setting $h_{j+1}:=f$ and $N'_{j+1}=h_{j+1}(M_{j+1})$.Ê Let $M'_{j+1}$ be
 an of $M'_j$ containing $N'_{j+1}$ of cardinality
 $\lambda$.Ê If possible, choose $M'_{j+1}$ to be universal over $M'_j$.
 ÊFix $\langle a_\zeta\mid \zeta\in S_{j+1}\rangle$ some
 enumeration of $\{a\in M'_{j+1}\mid a\models r\}$.

 In the event that $a_j\notin N'_j$, we need to be more careful with the
 amalgamation.Ê Let $f$ and $M^*$ be as in the diagram above.Ê Let us
 rewrite this diagram as
 \[
 \xymatrix{\ar @{} [dr] M_{j+1}
 \ar[r]^{f}Ê &M^{*} \\
 M_j \ar[u]^{\id} \ar[r]_{h_j}
 & N'_{j} \ar[u]_{\id} \ar[r]_{\id}
 & M'_j}
 \]
 Since $a_j\notin N'_j$, we have that $a_j\in
 M'_j\backslash N'_j$.Ê Thus $\tp(a_j/N'_j,M'_j)$ is non-algebraic.
 Now lets compare this to the realization $c_j$ of $r$ in
 $M_{j+1}\backslash M_j$. Notice that
 $f(\tp(c_j/M_j,M_{j+1}))=\tp(f(c_j)/N'_j,M^*)$ is non-algebraic.Ê Since
 $h_0=\id_{M_0}$, $\tp(f(c_j)/N'_j,M^*)$ is alsoÊ an extension of $r$.
 By the minimality of $r$ we can conclude that
 $$\tp(a_j/N'_j,M'_j)=\tp(f(c_j)/N'_j,M^*).$$
 So we can find an amalgam $M^{**}$ such that the following diagram
 commutes

 $$ \xymatrix{\ar @{} [dr] M_{j+1}
 \ar[r]^{f}Ê &M^*\ar[r]^{g} &M^{**} \\
 M_j \ar[u]^{\id} \ar[r]_{h_j}
 & N'_{j} \ar[u]_{\id} \ar[r]_{\id}&
 M'_j \ar[u]_{\id}}
 $$
 and $g(f(c_j))=a_j$.Ê Let $h_{j+1}:=g\circ f$ and $N'_{j+1}:=
 h_{j+1}(M_{j+1})$.Ê Let $M'_{j+1}$ be anÊ extension of $M'_j$
 cardinality $\lambda$ containing $N'_{j+1}$.Ê If possible, choose
 $M'_{j+1}$ to be universal over $M'_j$.Ê Fix
 $\langle a_\zeta\mid\zeta\in S_i\rangle$ some enumeration of $\{a\in
 M'_i\mid a\models r\}$.
 This completes the construction.

 Let $N':=\Union_{i<\alpha}N'_i$, $M':=\Union_{i<\alpha}M'_i$ and
 $h:=\Union_{i<\alpha}h_i$.Ê Notice that $M'$ realizes $q$
 and $h:M\cong N'$ with $h\restriction M_0=\id_{M_0}$.Ê We will show that
 $N'=M'$ in order to conclude that $M$ also realizes $q$.Ê Suppose not.
 Then
 $N'\prec_{\K}M'$ and we can fix $a\in M'\backslash N'$.Ê Since there are
 no $(r,\lambda)$-Vaughtian pairs, we can choose $a$ such that $a\models
 r$.

 By the definition of
 $M'$, there is an
 $i<\alpha$ such that
 $a\in M'_i$.Ê Then $a=a_\zeta$ for some $\zeta\in S_i$.Ê At stage,
 $\zeta+1$, we made sure that $a=a_\zeta\in N'_{\zeta+1}\subseteq N'$.
 This contradicts our choice of $a$.
 \end{proof}

 We now restate Theorem \ref{step up}
 \begin{theorem}
 Suppose that $\K$ has arbitrarily large models, is $\chi$-tame and
 satisfies the amalgamation property. If $\lambda\geq\chi\geq\LS(\K)$
 and
 $\K$ is categorical in both
 $\lambda$ and
 $\lambda^+$ then
 $\K$ is categorical in $\lambda^{++}$.
 \end{theorem}









 \begin{proof}[Proof of Theorem \ref{step up}]

 We will show that for every
 $N\in\K_{\lambda^{++}}$ and every $M\prec_{\K}N$ of cardinality
 $\lambda^+$, $N$ realizes every type over $M$.

 Let $M\prec_{\K}N$ have cardinality $\lambda^+$.
 First notice that Proposition \ref{locally min exist prop} and
 categoricity in
 $\lambda^+$ guarantees that there existsÊ a rooted minimal
 $r\in\gaS(M)$. By Corollary \ref{many
 realizations}, we know that $N$ realizes $r$ $\lambda^{++}$-times.

 Let $\alpha<\lambda^+$ be such that
 $\alpha=\lambda^+\cdot\alpha$. By the
 Downward-L\"{o}wenheim Skolem Axiom of AECs, we can construct a
 $\prec_{\K}$-increasing and continuous chain of models
 $\langle M_i\prec_{\K}N\mid i<\alpha\rangle $ such that $M=M_0$
 for every $i<\alpha$, we can fix $a_i\in M_{i+1}\backslash M_i$ realizing
 $r$.Ê This construction is possible since there are $\lambda^{++}$-many
 realizations of $r$ to choose from.Ê By Theorem \ref{many realizations of
 minimal imply sat},
 $\Union_{i<\alpha}M_i$ realizes every type over $M$.

 \end{proof}

 We now can derive an upward categoricity transfer theorem

 \begin{corollary}[Categoricity Transfer for Tame AECs]\label{cor categ in
 omega} Suppose that $\K$ has arbitrarily large models, satisfies the
 amalgamation property and is $\chi$-tame with
 $\chi\geq\LS(\K)$.
 Suppose that $\lambda\geq\max\{\chi,\LS(\K)^+\}$.
 If
 $\K$Ê is categorical in both $\lambda^+$ and $\lambda$, then
 $\K$ is categorical in every
 $\mu$ with $\lambda\leq\mu$.
 \end{corollary}

 categ in
 chain of
 %

 \begin{proof}[Proof of Corollary \ref{cor categ in
 omega}]
 Let $\alpha$ be such that $\lambda=\aleph_\alpha$.
 We will prove that $\K$ is categorical in $\aleph_\beta$Ê for all
 $\beta\geq\alpha+2$.Ê The base case is Theorem \ref{step up}.
 For $\beta
 =\gamma+2$,
 ÊTheorem \ref{step up} and the induction hypothesisÊ give us
 that
 $\K$ is categorical in $\aleph_{\gamma+2}$.

 \underline{Deriving categoricity in $\aleph_{\beta+1}$ where $\beta$ is a
 limit ordinal $>\alpha$}. Assume that
 $\K$ is categorical in every $\mu$ with
 $\lambda\leq\mu\leq\aleph_\beta$.Ê We need to show that $\K$ is
 categorical in $\aleph_{\beta+1}$.
 Let $N$ have cardinality $\aleph_{\beta+1}$.Ê Fix $M\prec_{\K}N$ with
 cardinality $\aleph_\beta$.Ê By Proposition \ref{locally min exist prop},
 there is a rooted minimal type $r\in\gaS(M)$.

 %
 %
 %

 Let $M'\in\K_\kappa$ with $\lambda\leq\kappa<\aleph_\beta$ be a root of
 $r$. Suppose that there is a $(r,\aleph_\beta)$-Vaughtian pair.Ê Then by
 Theorem \ref{vp transfer}, there is a $(r\restriction
 M',\kappa)$-Vaughtian pair.Ê Our induction hypothesis tells us that $\K$
 is categorical in $\kappa$ and
 $\kappa^+$.Ê Thus there are no $(r\restriction M',\kappa)$-Vaughtian
 pairs.Ê And we can conclude there are no $(r,\aleph_\beta)$-Vaughtian
 pairs.

 We now see that $N$ realizes
 $r$Ê $\aleph_{\beta+1}$-many times.Ê Since there are enough realizations
 of $r$ to go around, we can construct an increasing and continuous
 chain of models $\langle M_i\prec_{\K}N\mid i<\gamma\rangle$ of
 cardinality $\aleph_\beta$ such that for every
 $i<\gamma$, there is a $c_i\in M_{i+1}\backslash M_i$ realizing $r$
 andÊ $\gamma$ is a limit ordinal
 $<\aleph_\beta$ satisfying $\gamma=\aleph_\beta\cdot\gamma$.Ê Now by
 Theorem \ref{many realizations of minimal imply sat} we see that $N$ must
 realize every type over $M$.

 \underline{Deriving categoricity in $\aleph_\beta$ for $\beta$ a limit
 ordinal}.
 Assume that $\K$
 is categorical in all $\mu$ with $\lambda\leq\mu<\aleph_\beta$. For
 this case, it suffices to show that the every model of
 cardinality $\aleph_{\beta}$ is Galois-saturated.Ê ÊGiven
 $N\in\K_{\aleph_{\beta}}$ and
 $M\prec_{\K}N$ be a model of cardinality $\aleph_{\gamma}$ for
 someÊ $\gamma>\alpha$.Ê Ê Let
 $p\in
 \gaS(M)$ be given.Ê By the Downward L\"{o}wenheim Skolem axiom of AECs,
 we may find $N'\in\K_{\aleph_{\gamma+1}}$ such that
 $M\prec_{\K}N'\prec_{\K}N$.Ê By the induction hypothesis $N'$ is
 Galois-saturated and realizes $p$.Ê Thus $N$ realizes
 $p$.
 \end{proof}




 Combining Corollary \ref{cor categ in omega} and Fact \ref{statement of
 Sh394} yields

 \begin{corollary}\label{final corollary}

 Suppose $\K$ is a $\chi$-tame abstract elementary class satisfying the
 amalgamation and joint embeddingÊ properties.
 Let
 $\mu_0:=\Hanf(\K)$.Ê If $\chi\leq\beth_{(2^{\mu_0})^+}$ and
 $\K$ is categorical in some
 $\lambda^+>\beth_{(2^{\mu_0})^+}$, then $\K$ is categorical in $\mu$ for
 all
 $\mu>\beth_{(2^{\mu_0})^+}$.

 \end{corollary}

 %


induction

 \bigskip

 \section{Implications and Open Problems} \label{s:implications}

 The Hart-Shelah examples \cite{HaSh} (an alternative exposition is in
 chapter 19 of
 \cite{Ba1}) have
 arbitrary large models and are categorical in several successive
 cardinals but fail
 to be categorical in some larger cardinals.Ê By Theorem
 \ref{step up}Ê the AEC induced by $\phi_n$\footnote{Let $\phi_n$
 represent an
 $L_{\omega_1,\omega}$ formula that axiomatizes the Hart-Shelah example
 which is categorical in $\aleph_0,\dots,\aleph_n$ but not categorical in
 some larger cardinality.} is not
 $\chi$-tame for any
 $\chi<\aleph_n$ or fails to have the amalgamation property.


%
%
 Baldwin, David Kueker, Grossberg and VanDieren have begun extending the
 results for categorical, tame AECs to stable, tame AECs.Ê After the
 presentation of the results from this paper at the Bogot\'{a} Meeting in
 Model Theory 2003, Lessmann and Tapani Hytinnen have explored the
 implications of our arguments in more specific contexts.

 The following paragraph was added at the request of the referee:
 During the summer of 2004, Lessmann asked us if it was possible to prove
 the conclusion of Theorem
 \ref{step up} by assuming only categoricity in $\lambda^+$.
This startedÊ
a sequence of
 emails with him in which we discussed the issues and provided him
explanations of our ideas.Ê We were under the impression that
 this will be a collaboration between the three of us.Ê Unfortunately he did not
inform us
 that he did not intend to collaborate with us.Ê A few weeks after that
exchange, to our
 surprise, we received from him his paper where he proves an upward
categoricity result
 from a successor assuming that $\K$ has arbitrarily large models, has that
amalgamation
 property, $\LS(\K)=\aleph_0$ and $\aleph_0$-tameness  using ideas
 provided to him through our e-mail exchanges.

 Immediately after reciving his paper we have eliminated the assumptions that
 $\LS(\K)=\aleph_0$ and that $\K$ is $\aleph_0$-tame and proved:

 \begin{theorem}[\cite{GrVa2}] Suppose $\K$ is a $\chi$-tame AEC with the
amalgamation
 property with arbitrarily large models.
 If $\K$ is categorical in $\lambda^+$ for some $\lambda>(\LS(\K)+\chi)$ then
$\K$ is
 categorical in all
 $\mu>(\LS(\K)+\chi)$

 \end{theorem}


 \end{document}